\documentclass[12pt]{article}
\usepackage{graphicx}
\usepackage{amsfonts,amsmath}
\usepackage[mathscr]{eucal}
\usepackage{amssymb}
\usepackage{amsthm}
\usepackage{bbold}
\theoremstyle{plain}
\newtheorem{thm}{Theorem}

\newcommand{\be}{\begin{eqnarray}}
\newcommand{\ee}{\end{eqnarray}}
\newcommand{\bc}{\begin{center}}
\newcommand{\ec}{\end{center}}
\newcommand{\nn}{\nonumber \\}
\newcommand{\lb}{\label}
\newcommand{\pd}{\partial}
\newcommand{\p}[1]{(\ref{#1})}
\newcommand{\vecg}[1]{\mbox{\boldmath $#1$}}

\begin{document}

\begin{titlepage}

\vspace*{0.2cm}

\renewcommand{\thefootnote}{\star}
\begin{center}

{\LARGE\bf  A few comments on (hyper)k\"ahler geometry}

\vspace{2cm}

{\Large A.V. Smilga} \\

\vspace{0.5cm}

{\it SUBATECH, Universit\'e de
Nantes,  4 rue Alfred Kastler, BP 20722, Nantes  44307, France. }

\end{center}
\vspace{0.2cm} \vskip 0.6truecm \nopagebreak

   \begin{abstract}
\noindent  
In this note, we make two methodical observations. 

\begin{enumerate}
\item We prove in a simple explicit way that 
a necessary and sufficient condition for a K\"ahler manifold to be hyperk\"ahler is
 \be
\label{many-heaven}
h_{i\bar k} h_{j\bar l } \Omega^{\bar k \bar l} \ =\ C \Omega_{ij},
 \ee
where $h_{i\bar k}$ is a complex metric, $\Omega$ is a symplectic matrix and $C$ is a positive constant.
\item
The procedure of K\"ahler reduction includes two stages. On the first stage, a K\"ahler manifold of dimension $2n$ is reduced to a $(2n-1)$- dimensional manifold, while on the second stage, one arrives at a K\"ahler manifold of dimension $2(n-1)$. We note that this second stage has the meaning of Hamiltonian reduction. We illustrate the procedure by discussing a simple toy model when $\mathbb{R}^3 \times S^1$ is reduced down to $S^2$.

We elucidate also hyperk\"ahler reduction of  $\mathbb{R}^8$  down to the Taub-NUT metric.
\end{enumerate}

   \end{abstract}

\end{titlepage}

\section{Multidimensional heavenly equation}

Consider a K\"ahler manifold with the metric
\be h_{i\bar k} \ =\ \frac {\pd^2 {\cal K}(z, \bar z)}{\pd z^i \pd \bar z^{\bar k}} \ee

of complex dimension 2. We know since \cite{Pleb} that this manifold is hyperk\"ahler iff
the condition\footnote{For any $d \geq 2$, the condition \p{det-h} stipulates that the manifold is Ricci-flat (see e.g. Ref. \cite{Candelas}). Ricci flatness is a necessary condition for the manifold to be hyperk\"ahler. When $d=2$, it is also a sufficient condition.}

\be
\lb{det-h}
\det (h_{i\bar k}) \ =\ C >0
 \ee
is fulfilled. This is the famous {\it heavenly equation}. By rescaling the metric, we can set $C=1$.

To give one of many possible elementary proofs of this fact, note first that 
a Hermitian $2\times 2$ matrix with unit determinant can be presented as
 \be
\lb{exp-vsigma}
h_{i\bar k} \ =\ \left[ \exp(v_a \sigma_a) \right]_{j\bar k} = \cosh |\vecg{v}| \,\delta_{i\bar k} + \frac {\sinh |\vecg{v}|}{|\vecg{v}|}v_a (\sigma_a)_{i\bar k}.
 \ee
Here $\sigma_a$ are the Pauli matrices and  $v_{1,2,3}$ are real numbers. The validity of \p{det-h} can be checked quite directly.

Were $v_a$ purely imaginary, the matrix \p{exp-vsigma} would be an element of $SU(2)$.
As $v_a$ are real, $h_{j\bar k}$ belongs not to $SU(2)$, but  to a coset: the complexified version of $SU(2)$ [i.e. $SL(2, \mathbb{C})$] factorized over its compact form. 
 It is natural then that the holonomy group of such a manifold is $SU(2)$ rather than $U(2)$.
 
The proof of the last assertion is easy. If $\det h$ = 1, one can choose the vielbeins such as $\det e = \det \bar e = 1$. Consider the K\"ahler spin connection\footnote{We use the same  notation  as in Ref. \cite{glasses}. The lowercase Latin letters from the middle of the alphabet are the indices of world complex vectors $w^m$, the barred indices refer to complex conjugated vectors $\bar w^{\bar m}$. The lowercase letters from the beginning of the aphabet are the tangent space complex indices. The capital Latin letters are used for the world and tangent space indices of real vectors.}
\be
              \lb{omega-Kahl}
              \omega_{\bar b a, m}  = -  \omega_{a \bar b, m} \ =\ e^{\bar n}_{\bar b} \partial_m e^{\bar a}_{\bar n}\, , \quad
               \omega_{ b \bar a, \bar m}  = - \omega_{\bar a b, \bar m}\ =\ e^{n}_{b} \partial_{\bar m} e^{a}_{n}. \quad
              \ee 
We see that
\be
\omega_{a\bar a, \bar m}  = - \omega_{\bar a a, \bar m}  = e^n_a \pd_{\bar m} e^a_n = \pd_{\bar m} (\ln \det e) \ =\ 0\,.
\ee
Similarly, $\omega_{a\bar a,  m} = 0 $. Then $\omega_{AB, M}dx^M \in su(2)$, the same concerns the curvature form
\be
R_{AB} \ =\ d \omega_{AB} + \omega_{AC} \wedge \omega_{CB},
\ee
and the holonomy group of the manifold is $SU(2)$ rather that $U(2)$ what would be the case for a generic K\"ahler 4-dimensional  manifold. It follows\footnote{This is a well-known mathematical fact, see e.g. Theorem 3.2 in Ref. \cite{glasses}.} that the manifold is hyperk\"ahler.

We go now to the multidimensional case and prove the theorem:
\begin{thm}
A K\"ahler manifold with the metric $h_{i\bar k}$ is hyperk\"ahler if the  relation \p{many-heaven} holds.
 $\Omega_{ij}$ is  the symplectic matrix  that can be chosen in the block-diagonal form $\Omega = {\rm diag} (\varepsilon, \ldots, \varepsilon)$ with 
$$ \varepsilon \ =\ \left(\begin{array}{cc}0 & 1 \\ -1 & 0
 \end{array}   \right).$$
\end{thm}

The relation \p{many-heaven} is a natural generalization of the heavenly equation \p{det-h} for the multidimensional case. 
It has the same form as the condition $U \Omega U^T = \Omega$, which defines an $Sp(n)$ matrix $U$, but the Hermitian metric $h_{j\bar k}$ is not an element of $Sp(n)$, but rather, in a full analogy with \p{exp-vsigma}, has the form
\be
\lb{exp-vt}
h_{j\bar k} \ =\ \left[ \exp(v_a t_a) \right]_{j\bar k}\,,
 \ee
where $t_a$ are the Hermitian generators of $Sp(n)$ satisfying the condition $(t_a)^T \Omega + \Omega t_a = 0$ and $v_a$ are real.
Thus, $h_{j\bar k}$ is an element of the coset: the complexified version of $Sp(n)$ factorized over its compact form. It is natural then that the holonomy group of the manifold is $Sp(n)$  and hence the manifold is hyperk\"ahler. 

This theorem was proved earlier in Refs. \cite{Bridgeland,Dunajski}, but we suggest an elementary explicit proof. 

\begin{proof}
We choose the complex structures having the following nonzero components:
\be
&&I_{m \bar n} = -i h_{m\bar n} \ \Longrightarrow  I_m^n = -i \delta_m^n, \ I_{\bar m}^{\bar n} = i \delta_{\bar m}^{\bar n}, \nn
 && J_{mn} = \Omega_{mn}, \ J_{{\bar m} {\bar n} } = \Omega_{{\bar m} {\bar n}},
\quad K_{mn} = -i\Omega_{mn}, \ K_{{\bar m} {\bar n} } =  i\Omega_{{\bar m} {\bar n}}.
\ee

Using the condition \p{many-heaven}, it is easy to show that they are quaternionic:
\be 
I^2 = J^2 = K^2 = -\mathbb{1}, \quad IJ = K, \, JK = I, \, KI = J.
 \ee
To prove the theorem, one should also show that they are  covariantly constant. This is obvious for   the conventional K\"ahler complex structure $I$. 

It is convenient to introduce the combinations $J\pm iK$. 
$J+iK$ has only holomorphic and $J - iK$ only antiholomorphic components. 
The ordinary derivatives of $J \pm iK$ vanish, and we have to prove that
\be \nabla_P (J\pm iK)_{MN} = \Gamma^S_{PM} (J \pm  iK)_{SN} + \Gamma^S_{PN} (J \pm iK)_{MS} = 0.
\ee

It is well known that the only nonzero Christoffel symbols for K\"ahler manifolds are purely holomorphic, $\Gamma^s_{mn} = h^{\bar q s} \partial_m h_{n \bar q}$, or antiholomorphic.  Hence  $ \nabla_{\bar p} (J+iK)_{mn} =2 \partial_{\bar p} \Omega_{mn} = 0$. Now,
\be \nabla_p (J+iK)_{mn} \ = 2 \ h^{\bar q s} \partial_p h_{m \bar q} \Omega_{sn} - (m \leftrightarrow n).
\ee
The matrix  $(X_p)_m^s =  2h^{\bar q s} \partial_p h_{m \bar q} $ is an element of $sp(n)$, as follows from \p{exp-vt}. Then,

\be && \nabla_p (J+iK)_{mn}  \ = \   (X_p)_m^s \Omega_{sn} - (X_p)_n^s \Omega_{sm} \nn
&& = \ 
[X_p \Omega + \Omega (X_p)^T]_{mn}  \ = \ 0.
\ee

The property $\nabla_P (J-iK)_{MN} = 0$ is proved in the same way.

\end{proof}

This reasoning can be applied also in the four-dimensional case to give an alternative derivation of the heavenly equation 
\p{det-h}.

\section{K\"ahler reduction: a toy model}
The method of K\"ahler reduction is based on the observations made half-a-century ago in Ref. \cite{Marsden}. It was developped and generalized for hyperk\"ahler manifolds in Ref. \cite{Hitchin} and further in Ref. \cite{Rychenkova}. This method can be applied when the parent  metric enjoys certain isometries. It is especially practical in the hyperk\"ahler case, allowing one to derive many nontrivial hyperk\"ahler metrics. This procedure is not so transparent, however, and we found it useful to see how it works in simple cases.

As probably the simplest nontrivial example, consider the space $\mathbb{R}^2 \times (\mathbb{R} \times S^1)$ with the flat metric expressed in the form
 \be
\lb{metr-R4}
     ds^2 \ =\ dr^2 + r^2 d\phi^2 + dx^2 + a^2 d\theta^2\,,
\ee
where $\phi, \theta \in (0,2\pi)$ and $a$ is the radius of the circle. One of the isometries of this manifold is a 
simultaneous shift of $\phi$ and $\theta$:
  \be
\lb{alpha}
\phi \to \phi + \alpha, \qquad \theta \to \theta + \alpha \,.
 \ee

The difference
\be
\lb{chi}
\chi \ =\ \phi - \theta
 \ee
stays invariant under \p{alpha}, and the same concerns the K\"ahler form
 \be
\lb{om}
\omega = \omega_{MN} dx^M \wedge dx^N \ =\ r dr \wedge d\phi + a dx \wedge d\theta\,.
 \ee
Consider the 1-form $\omega_\alpha$ obtained by a {\it contraction} of the form \p{om} with the vector $V^M = (0_r, 1_\phi, 0_x, 1_\theta)$ generating the transformation \p{alpha}:
 \be
\lb{om-alpha}
\omega_\alpha \ =\ \omega_{MN} dx^M V_N \ =\ r dr + a dx\,.
 \ee
The form \p{om-alpha} is obviously closed and also exact, $\omega_\alpha = d\mu_\alpha$, with
 \be
\lb{mu-alpha}
 \mu_\alpha \ =\ \frac {r^2}2 + ax \,.
 \ee
The function $\mu_\alpha$ is called a {\it moment map}.

To perform K\"ahler reduction, we set $\mu_\alpha = 0$ and substitute
\be
\lb{x-via-r}
x \ =\ - \frac {r^2}{2a}
\ee 
into \p{metr-R4}. The metric acquires the form
\be
\lb{3D-toy}
ds^2 \ =\ dr^2 \left( 1 + \frac{r^2}{a^2} \right) + r^2 (d\theta + d\chi)^2 + a^2 d\theta^2 \,.
 \ee
Next, we factorize the three-dimensional manifold thus obtained with respect to the isometry \p{alpha}. To this end, we represent the metric as
\be
\lb{3D-toyy}
\!\!\!\!\!\!\!\!\! ds^2 = dr^2 \left( 1 + \frac{r^2}{a^2} \right)  + \frac{a^2r^2}{r^2 + a^2} d\chi^2 + (r^2 + a^2)\left(d\theta + \frac{r^2}{r^2 + a^2} d\chi \right)^2
\ee
and {\it suppress} the last term with complete square.
We are left with 
  \be
\lb{red-metr}
ds^2 \ =\ dr^2 \left( 1 + \frac{r^2}{a^2} \right)  + \frac{a^2r^2}{r^2 + a^2} d\chi^2 \,.
  \ee

The last operation has a dynamical interpretation. Consider the Lagrangian describing the motion of a particle along the manifold  \p{3D-toy}. It reads
\be
L \ =\ \left(1 + \frac {r^2}{a^2} \right) \frac {\dot r^2}2 + \frac {r^2}2 \dot \chi^2  + \frac {r^2 + a^2}2 \dot \theta^2
+ r^2 \dot \chi \dot \theta\,.
\ee
The corresponding canonical Hamiltonian is
\be
H \ =\ \frac {p_r^2}{2(1 + r^2/a^2)} + \frac {p_\chi^2 (1 + r^2/a^2)}{2r^2} + \frac {p_\theta^2}{2a^2} - \frac {p_\chi p_\theta}{a^2}\,.
\ee
The Poisson bracket $\{p_\theta, H\}_P$ vanishes, and we can impose the constraint $p_\theta = 0$. 
This gives the reduced Hamiltonian
 \be
H \ =\ \frac {p_r^2}{2(1 + r^2/a^2)} + \frac {p_\chi^2 (1+ r^2/a^2)}{2 r^2} \,.
 \ee
Its Legendre transform,
\be
L \ =\ \frac{\dot r^2}{2} \left( 1 + \frac {r^2}{a^2} \right) +  \frac{\dot \chi^2}{2} \frac{r^2}{1 + r^2/a^2}\,, 
\ee  
corresponds to the metric \p{red-metr}.

 This part of the procedure is nothing but  {\it Hamiltonian reduction} --- the method invented  by Dirac almost 100 years ago \cite{Dirac} and well known to physicists.\footnote{We do not touch here, however, upon the quantum aspect of the problem.}

The metric \p{red-metr} represents a  hemisphere, recembling a tea cup. Its Gaussian curvature is positive,
\be
\mathbb{R} \ =\ \frac {4a^4}{(r^2 + a^2)^3},
\ee
and the Gauss-Bonnet integral is
\be \chi \ = \ \frac 1{2\pi} \int d^2x \sqrt{g} \,\mathbb{R}  \ =\  4a^4\int_0^\infty 
\frac {rdr}{(r^2 + a^2)^3} \ =\ 1. \ee

  As any 2D manifold, this manifold is K\"ahler. Its K\"ahler form,
\be
\lb{om-tilde}
\tilde \omega = rdr \wedge d\chi\,,
 \ee
 {\it coincides} with \p{om}, bearing in mind \p{chi} and \p{x-via-r}.

In this case, it is not difficult  to see by
a direct calculation that the complex structure $I_{r\chi} = - I_{\chi r} = r$ associated with $\tilde \omega$ is covariantly constant.\footnote{It is true for any metric
$$ ds^2 \ =\ f(r) dr^2 + \frac {r^2}{f(r)} d\chi^2\,. $$}
 Moreover, mathematicians can prove (see  Ref. \cite{Marsden} and Theorem 3.1 in Ref. \cite{Hitchin}) that it is 
{\it always} so and that the complex structure on the quotient derived by applying the procedure described above to any K\"ahler manifold is covariantly constant so that the quotient is also K\"ahler.

\section{Hyperk\"ahler reduction: Taub-NUT}

Hyperk\"ahler reduction is K\"ahler reduction with respect to  all the three complex structures. It is possible to perform for hyperk\"ahler manifolds possessing an isometry leaving intact three K\"ahler forms. The simplest nontrivial example \cite{Rychenkova} is $\mathbb{R}^8$ or rather $\mathbb{R}^7 \times S^1  = \mathbb{R}^4 \times (\mathbb{R}^3 \times S^1)$ with the Taub--NUT manifold \cite{TAUB} as a quotient.\footnote{See also Ref. \cite{Gaeta} and Section 2.5 of Ref. \cite{Wit-GRG} for  pedagogical expositions.}

The flat $\mathbb{R}^4$ metric,
\be 
\lb{flat-R4-Cart}
ds^2 = dy_1^2 + dy_2^2 + dy_3^2 + dy_4^2,
 \ee
 admits three quaternionic complex structures. The corresponding K\"ahler forms are
\be
\lb{om-IJK}
\omega_I  \ = \ dy_1 \wedge dy_2 + dy_3 \wedge dy_4\,, \nn
 \omega_J \ = \ dy_1 \wedge dy_3 - dy_2 \wedge dy_4\,, \nn
\omega_K  \  = \ dy_1 \wedge dy_4 + dy_2 \wedge dy_3\,.
 \ee
These forms are invariant under the isometry
\be
\lb{alpha4}
&&y_1' = y_1 \cos  \alpha  - y_2 \sin  \alpha , \qquad y_2' = y_1 \sin  \alpha  + y_2 \cos \alpha \,, \nn
&&y_3' = y_3 \cos  \alpha  + y_4 \sin \alpha , \qquad y_4' =  -y_3 \sin \alpha  + y_4 \cos  \alpha \,.
\ee
Let us make the following variable change:
\be
\lb{var-change}
&x_1 \ =\ 2(y_1 y_4 + y_2 y_3)    , \qquad
x_2 \ =\ 2(  y_2 y_4 - y_1 y_3 ) , \nn
&x_3 \ =\ y_1^2 + y_2^2 - y_3^2 - y_4^2, \qquad
\Psi \ =\  -2\arctan \frac{y_1}{y_2},
 \ee
 The coordinates $\{\vecg{x}, \Psi\}$ play on $\mathbb{R}^4$ the same role as the polar coordinates on $\mathbb{R}^2$. In these terms, the metric \p{flat-R4-Cart} acquires the form \cite{GH-PLB}:
\be
\lb{flat-R4-Psi}
ds^2 \ =\ \frac r4 (d\Psi + \vecg{A} \cdot d\vecg{x} )^2 \ +\ \frac 1{4r} d\vecg{x} \cdot d\vecg{x}, 
\ee
where $r = |\vecg{x}|$ and 
\be
\lb{A-monop}
A_1 = \frac {x_2}{r(r + x_3)}, \qquad A_2 = - \frac {x_1}{r(r + x_3)}, \qquad A_3  =0 
\ee
is nothing but the field of the {\it Dirac monopole} \cite{Dirac-monop} written in a particular gauge.\footnote{One of the ways to derive this nontrivial fact is using a quaternionic representation of the metric --- see Eq.(11) in Ref. \cite{Rychenkova}.}

Then the isometry \p{alpha4} boils down to the shift  $\Psi \to \Psi + 2\alpha$, while the variables $x_j$ do not change.

Consider now the manifold $\mathbb{R}^4 \times (\mathbb{R}^3 \times S^1)$ choosing the metric \p{flat-R4-Psi} for the first factor and the Cartesian metric 
$ds^2 = d\vecg{X}^2 + a^2 d\theta^2$ for $\mathbb{R}^3 \times S^1$, $a$ being the radius of the circle. 
Consider the isometry\footnote{$\theta \in (0,2\pi)$ and $\Psi \in (0,4\pi)$.}
\be
\lb{iso-G}
G: \qquad \Psi \to \Psi + 2\alpha, \ \theta \to \theta + \alpha.
\ee
We choose the hyperk\"ahler triple in the form\footnote{A little ``twist" between the sectors $\mathbb{R}^4$ and $\mathbb{R}^3 \times S^1$ is a pure convention assuring the positive signs in the RHS of Eq. \p{om-dmu}.}  
\be
\lb{triple-R8}
\omega_I \ =\ dy_1 \wedge dy_2 + dy_3 \wedge dy_4 \ +\ dX_1 \wedge dX_2 + dX_3 \wedge (a d\theta)\,,\nn
\omega_J \ =\ dy_1 \wedge dy_3 - dy_2 \wedge dy_4 \ +\ dX_2 \wedge dX_3  + dX_1 \wedge (a d\theta)    \,, \nn
\omega_K \ =\ dy_1 \wedge dy_4 + dy_2 \wedge dy_3 \ -\  dX_1 \wedge dX_3  +  dX_2 \wedge (a d\theta)\,.
\ee

The forms \p{triple-R8} 
 are invariant under $G$.   The combination $\chi  = \Psi - 2\theta$ also stays invariant. 
At the next step, we determine the contractions of these forms with the vector field
\be
\lb{VM}
V^M \ =\ (-y_2, y_1, y_4, -y_3, \vecg{0}, 1)
\ee
defining the infinitesimal rotation \p{alpha4}, \p{iso-G}.  We obtain:
\be
\lb{om-dmu}
&&\omega^\alpha_I \ =\ y_1 dy_1 + y_2 dy_2 - y_3 dy_3 - y_4 dy_4  +a dX_3 \ =\ d\left(\frac{x_3}2 + aX_3\right)\,, \nn
&&\omega^\alpha_J \ =\ y_4 dy_1 + y_3 dy_2 + y_2 dy_3  + y_1 dy_4 + a dX_1 \ =\  d\left(\frac{x_1}2 + aX_1\right) \,, \nn
&&\omega^\alpha_K \ =\ - y_3 dy_1 + y_4 dy_2 - y_1 dy_3 + y_2 dy_4 + a dX_2 \ =\ d\left(\frac {x_2}2 + aX_2 \right).
\ee

In other words, these forms represent exterior derivatives of  the moment maps: 
\be
\lb{mom-maps}
\mu_I = \frac{x_3}2 + aX_3, \quad \mu_J = \frac{x_1}2 + aX_1 , \quad  \mu_K =  \frac{x_2}2 +   aX_2 \,.
\ee
To perform the hyperk\"ahler reduction, we impose on the flat $\mathbb{R}^7 \times S^1$ metric the constraints $\mu_{I,J,K} = 0$. This gives us a five-dimensional manifold with the metric 
 \be
\lb{metr-D5}
ds^2 \ =\ \frac 14 \left(\frac 1r + \frac 1{a^2}  \right) d\vecg{x}^2 + \frac r4 \big(d\chi + 2d\theta + \vecg{A} \cdot d\vecg{x} \big)^2 + a^2 d\theta^2\,,
\ee
where we also plugged $\Psi = \chi + 2\theta$. For a further factorization with respect to  the action of \p{iso-G}, we rewrite \p{metr-D5} distinguishing a complete square,
 \be 
ds^2 \ =\ \frac 14 \left(\frac 1r + \frac 1{a^2}  \right) d\vecg{x}^2 + \frac {\big(d\chi +  \vecg{A} \cdot d\vecg{x} \big)^2}{4\left(\frac 1r + \frac 1{a^2}\right)} \nn
+ \ (r + a^2) \left[ d\theta + \frac r{2(r + a^2)} \big(d\chi +  \vecg{A} \cdot d\vecg{x} \big) \right]^2\,,
 \ee
and cross the complete square out. What is left coincides with the metric of the Taub-NUT manifold expressed in the form \cite{Hawking}.

The hyperk\"ahler triple on the quotient may be derived starting from \p{triple-R8} in the same way as we derived \p{om-tilde} starting from \p{om} in the toy two-dimensional model. An explicit calculation was done in Ref. \cite{Gaeta}, which we follow. 
 At the first step, one should express the flat $\mathbb{R}^4$ triple \p{om-IJK} in terms of $\vecg{x}$ and $\Psi$, using \p{var-change}. After some calculation, one arrives at rather simple expressions: 
\be
\lb{R4-xPsi}
\omega_I^{\rm flat} \ &=&\ \frac 1{4r} dx_1 \wedge dx_2  + \frac 14 dx_3 \wedge (d\Psi + A_1 dx_1 + A_2 dx_2), \nn
\omega_J^{\rm flat} \ &=&\ \frac 1{4r} dx_2 \wedge dx_3  + \frac 14 dx_1 \wedge (d\Psi  + A_2 dx_2), \nn
\omega_K^{\rm flat} \ &=&\ \frac 1{4r} dx_3 \wedge dx_1  + \frac 14 dx_2 \wedge (d\Psi + A_1 dx_1)\,.
 \ee
where $A_{1,2}$ are the monopole vector potentials in \p{A-monop}.

Now, we substitute it in \p{triple-R8}. After setting there $X_j = -x_j/(2a)$ and $d\Psi = d\chi + 2d\theta$, one may observe that the terms $\propto dx_j \wedge d\theta$ cancel, and we arrive at a surprisingly simple answer:
  \be
\lb{IJK-TN}
\omega_I^{\rm TN} \ =&\ \frac 14 dx_1 \wedge dx_2 \left( \frac 1r + \frac 1{a^2} \right) + \frac 14 dx_3 \wedge (d\chi + \vecg{A} \cdot d\vecg{x}), \nn
\omega_J^{\rm TN} \ =&\ \frac 14 dx_2 \wedge dx_3 \left( \frac 1r + \frac 1{a^2} \right) + \frac 14 dx_1 \wedge (d\chi + \vecg{A} \cdot d\vecg{x}), \nn
\omega_K^{\rm TN} \ =&\ \frac 14 dx_3 \wedge dx_1 \left( \frac 1r + \frac 1{a^2} \right) + \frac 14 dx_2 \wedge (d\chi + \vecg{A} \cdot d\vecg{x}) \,.
\ee
The only difference compared to the flat structures \p{R4-xPsi} is the change $1/r \to 1/r + 1/a^2$.

We mention here also the paper \cite{Kelekci} where the Taub--NUT K\"ahler forms   were derived based on the expression for the Taub-NUT metric presented in the original papers \cite{TAUB} in terms of Cartan-Maurer forms. These formulas (see Eq. (3.10) in Ref. \cite{Kelekci}) are rather bulky, however. Probably, they reduce to \p{IJK-TN} after an appropriate variable change.

\vspace{1mm}

I am indebted to M. Dunajski, O. Kelekci, K. Krasnov and V. Rubtsov for illuminating discussions.

\end{document}